\documentclass{amsart}

\newtheorem{theorem}{Theorem}[section]

\theoremstyle{definition}

\theoremstyle{remark}

\numberwithin{equation}{section}

\begin{document}
\title[Identities of symmetry for generalized Bernoulli polynomials]{$\begin{array}{c}
         \text{Identities of symmetry for generalized Bernoulli polynomials}\\
       \end{array}$}

\author{dae san kim}
\address{Department of Mathematics, Sogang University, Seoul 121-742, Korea}
\email{dskim@sogong.ac.kr}
\thanks{}

\subjclass[]{}

\date{}

\dedicatory{}


\begin{abstract}
In this paper, we derive eight basic identities of symmetry in three
variables related to generalized Bernoulli polynomials and
generalized power sums. All of these are new, since there have been
results only about identities of symmetry in two variables. The
derivations of identities are based on the $p$-adic integral
expression of the generating function for the generalized Bernoulli
polynomials and the quotient of $p$-adic integrals that can be
expressed as the exponential generating function for the generalized
power sums.

\end{abstract}

\maketitle

\section{Introduction and preliminaries}
  Let $p$ be a fixed prime. Throughout this paper,
$\mathbb Z_{p},\mathbb Q_{p},\mathbb C_{p}$ will respectively denote
the ring of $p$-adic integers, the field of $p$-adic rational
numbers and the completion of the algebraic closure of $\mathbb
Q_{p}$. Let $d$ be a fixed positive integer. Then we let

\begin{align*}
X=X_{d}=\lim_{\overleftarrow{N}}\mathbb Z/dp^{N}\mathbb Z,
\end{align*}
and let $\pi: X\rightarrow \mathbb Z_{p}$ be the map given by the
inverse limit of the natural maps

\begin{align*}
\mathbb Z/dp^{N}\mathbb Z\rightarrow \mathbb Z/p^{N}\mathbb Z.
\end{align*}

If $g$ is a function on $\mathbb Z_{p}$, we will use the same
notation to denote the function $g\circ \pi$. Let $\chi:(\mathbb
Z/d\mathbb Z)^{*}\rightarrow\overline{\mathbb Q}^{*}$ be a
(primitive) Dirichlet character of conductor $d$.\\
Then it will be pulled back to $X$ via the natural map $X\rightarrow
\mathbb Z/d\mathbb Z$. Here we fix, once and for all, an imbedding
$\overline{\mathbb Q}\rightarrow \mathbb{C}_{p}$, so that $\chi$ is
regarded as a map of $X$ to $\mathbb{C}_{p}$(cf. \cite{K1}).

For a uniformly differentiable function $f:X\rightarrow\mathbb
{C}_{p}$, the $p$-adic integral of $f$ is defined(cf. \cite{T1}) by

\begin{align*}
\int_{X}f(z)d\mu(z)=\lim_{N\rightarrow\infty}\frac{1}{dP^{N}}\sum_{j=0}^{dp^{N}-1}f(j).
\end{align*}
Then it is easy to see that

\begin{align}\label{a}
\int_{X}f(z+1)d\mu(z)=\int_{X}f(z)d\mu(z)+f'(0).
\end{align}
\\
More generally, we deduce from (\ref{a}) that, for any positive
integer $n$,
\begin{align}\label{b}
\int_{X}f(z+n)d\mu(z)=\int_{X}f(z)d\mu(z)+\sum_{a=0}^{n-1}f'(a).
\end{align}

Let $|~|_{p}$ be the normalized absolute value of $\mathbb{C}_{p}$
such that $|p|_{p}=\frac{1}{p}$, and let

\begin{align}\label{c}
E=\{t\in \mathbb C_{p}\mid \mid t \mid_{p}< p^{- \frac{1}{p-1}}\}.
\end{align}
\\
Then, for each fixed $t\in E$, the function $e^{zt}$ is analytic on
$\mathbb{Z}_{p}$ and hence considered as a function on $X$, and, by
applying $(\ref{b})$ to $f$ with $f(z)=\chi(z)e^{zt}$, we get the
$p$-adic integral expression of the generating function for the
generalized Bernoulli numbers $B_{n,\chi}$ attached to $\chi$:
\begin{align}\label{d}
\int_{X}\chi(z)e^{zt}d\mu(z)=\frac{t}{e^{dt}-1}\sum_{a=0}^{d-1}\chi(a)e^{at}=\sum_{n=0}^{\infty}B_{n,\chi}\frac{t^{n}}{n!}~(t\in
E).
\end{align}
So we have the following $p$-adic integral expression of the
generating function for the generalized Bernoulli polynomials
$B_{n,\chi}(x)$ attached to $\chi$:
\begin{align}\label{e}
\int_{X}\chi(z)e^{(x+z)t}d\mu(z)=\frac{te^{xt}}{e^{dt}-1}\sum_{a=0}^{d-1}\chi(a)e^{at}=\sum_{n=0}^{\infty}B_{n,\chi}(x)\frac{t^{n}}{n!}~(t\in
E,~x\in \mathbb Z_{p}).
\end{align}
Also, from $(\ref{a})$ we have:
\begin{align}\label{f}
\int_{X}e^{zt}d\mu(z)=\frac{t}{e^{t}-1}.
\end{align}
Let $S_{k}(n,\chi)$ denote the $k$th generalized power sum of the
first $n+1$ nonnegative integers attached to $\chi$, namely
\begin{align}\label{g}
S_{k}(n,\chi)=\sum_{a=0}^{n}\chi(a)a^{k}=\chi(0)0^{k}+\chi(1)1^{k}+\cdots+\chi(n)n^{k}.
\end{align}
From $(\ref{d})$, $(\ref{f})$, and $(\ref{g})$, one easily derives
the following identities: for $w\in \mathbb Z_{>0}$,
\begin{align}
\label{h}
\frac{wd\int_{X}\chi(x)e^{xt}d\mu(x)}{\int_{X}e^{wdyt}d\mu(y)}=&\frac{e^{wdt}-1}{e^{dt}-1}\sum_{a=0}^{d-1}\chi(a)e^{at}\\
\label{i}
&=\sum_{a=0}^{wd-1}\chi(a)e^{at}\\
\label{j}
&=\sum_{k=0}^{\infty}S_{k}(wd-1,\chi)\frac{t^{k}}{k!}~(t\in E).
\end{align}
In what follows, we will always assume that the $p$-adic integrals
of the various (twisted) exponential functions on $X$ are defined
for $t\in E$ (cf. $(\ref{c})$), and therefore it will not be
mentioned.

\cite{DR1}, \cite{H1}, \cite{T2}, \cite{Tu1} and \cite{Y1} are some
of the previous works on identities of symmetry involving Bernoulli
polynomials and power sums. For the brief history, one is referred
to those papers.

In this paper, we will produce 8 basic identities of symmetry in
three variables $w_{1},w_{2},w_{3}$ related to generalized Bernoulli
polynomials and generalized power sums (cf. (\ref{l1})-(\ref{n1}),
(\ref{p1}), (\ref{r1})-(\ref{u1})). All of these seem to be new,
since there have been results only about identities of symmetry in
two variables in the literature(\cite{T3}). The following is stated
as Theorem 2 and  an example of the full six symmetries in $w_{1},
w_{2},w_{3}$.

\begin{equation*}
\begin{split}
&\sum_{k+l+m=n}\binom{n}{k,l,m}B_{k,\chi}(w_{1}y_{1})B_{l,\chi}(w_{2}y_{2})S_{m}(w_{3}d-1,\chi)w_{1}^{l+m}w_{2}^{k+m}w_{3}^{k+l-1}\\
=&\sum_{k+l+m=n}\binom{n}{k,l,m}B_{k,\chi}(w_{1}y_{1})B_{l,\chi}(w_{3}y_{2})S_{m}(w_{2}d-1,\chi)w_{1}^{l+m}w_{3}^{k+m}w_{2}^{k+l-1}\\
=&\sum_{k+l+m=n}\binom{n}{k,l,m}B_{k,\chi}(w_{2}y_{1})B_{l,\chi}(w_{1}y_{2})S_{m}(w_{3}d-1,\chi)w_{2}^{l+m}w_{1}^{k+m}w_{3}^{k+l-1}\\
=&\sum_{k+l+m=n}\binom{n}{k,l,m}B_{k,\chi}(w_{2}y_{1})B_{l,\chi}(w_{3}y_{2})S_{m}(w_{1}d-1,\chi)w_{2}^{l+m}w_{3}^{k+m}w_{1}^{k+l-1}\\
=&\sum_{k+l+m=n}\binom{n}{k,l,m}B_{k,\chi}(w_{3}y_{1})B_{l,\chi}(w_{2}y_{2})S_{m}(w_{1}d-1,\chi)w_{3}^{l+m}w_{2}^{k+m}w_{1}^{k+l-1}\\
=&\sum_{k+l+m=n}\binom{n}{k,l,m}B_{k,\chi}(w_{3}y_{1})B_{l,\chi}(w_{1}y_{2})S_{m}(w_{2}d-1,\chi)w_{3}^{l+m}w_{1}^{k+m}w_{2}^{k+l-1}.\\
\end{split}
\end{equation*}

The derivations of identities are based on the $p$-adic integral
expression of the generating function for the generalized Bernoulli
polynomials in (\ref{e}) and the quotient of integrals in
(\ref{h})-(\ref{j}) that can be expressed as the exponential
generating function for the generalized power sums. These abundance
of symmetries would not be unearthed if such $p$-adic integral
representations had not been available. We indebted this idea to the
paper \cite{T3}.

\section{Several types of quotients of $p$-adic integrals}

Here we will introduce several types of quotients of $p$-adic
integrals on $X$ or $X^{3}$ from which some interesting identities
follow owing to the built-in symmetries in $w_{1}, w_{2}, w_{3}$. In
the following, $w_{1}, w_{2}, w_{3}$ are positive integers and all
of the explicit expressions of integrals in (\ref{l}), (\ref{n}),
(\ref{p}), and (\ref{r}) are obtained from the identities in
(\ref{d}) and (\ref{f}). To ease notations, from now on we suppress
$\mu$ and denote, for example, $d\mu(x)$ simply by $dx$.
\\
\\
(a) Type $\Lambda_{23}^{i}$ (for $i=0,1,2,3$)\\
\begin{flalign*}
I(\Lambda_{23}^{i})
\qquad\qquad\qquad\qquad\qquad\qquad\qquad\qquad\qquad\qquad\qquad\qquad\qquad\qquad\qquad\qquad
\end{flalign*}
\begin{align}
\label{k}
&=\frac{d^{i}\int_{X^{3}}\chi(x_{1})\chi(x_{2})\chi(x_{3})e^{(w_{2}w_{3}x_{1}+w_{1}w_{3}x_{2}+w_{1}w_{2}x_{3}+w_{1}w_{2}w_{3}(\sum_{j=1}^{3-i}y_{j}))t}dx_{1}dx_{2}dx_{3}}{(\int_{X}e^{dw_{1}w_{2}w_{3}x_{4}t}dx_{4})^{i}}\\
\label{l}
&=\frac{(w_{1}w_{2}w_{3})^{2-i}t^{3-i}e^{w_{1}w_{2}w_{3}(\sum_{j=1}^{3-i}y_{j})t}(e^{dw_{1}w_{2}w_{3}t}-1)^{i}}{(e^{dw_{2}w_{3}t}-1)(e^{dw_{1}w_{3}t}-1)(e^{dw_{1}w_{2}t}-1)}
\end{align}
\begin{align*}
\qquad \qquad \qquad \qquad \times (\sum_{a=0}^{d-1}\chi(a)e^{aw_{2}w_{3}t})(\sum_{a=0}^{d-1}\chi(a)e^{aw_{1}w_{3}t})(\sum_{a=0}^{d-1}\chi(a)e^{aw_{1}w_{2}t})
\end{align*}
\\
\\
(b) Type $\Lambda_{13}^{i}$ (for $i=0,1,2,3$)
\begin{align}
\label{m}
I(\Lambda_{13}^{i})&=\frac{d^{i}\int_{X^{3}}\chi(x_{1})\chi(x_{2})\chi(x_{3})e^{(w_{1}x_{1}+w_{2}x_{2}+w_{3}x_{3}+w_{1}w_{2}w_{3}(\sum_{j=1}^{3-i}y_{j}))t}dx_{1}dx_{2}dx_{3}}{(\int_{X}e^{dw_{1}w_{2}w_{3}x_{4}t}dx_{4})^{i}}\\
\label{n}
&=\frac{(w_{1}w_{2}w_{3})^{1-i}t^{3-i}e^{w_{1}w_{2}w_{3}(\sum_{j=1}^{3-i}y_{j})t}(e^{dw_{1}w_{2}w_{3}t}-1)^{i}}{(e^{dw_{1}t}-1)(e^{dw_{2}t}-1)(e^{dw_{3}t}-1)}
\end{align}
\begin{align*}
\qquad \qquad \qquad \qquad\qquad \times (\sum_{a=0}^{d-1}\chi(a)e^{aw_{1}t})(\sum_{a=0}^{d-1}\chi(a)e^{aw_{2}t})(\sum_{a=0}^{d-1}\chi(a)e^{aw_{3}t})
\end{align*}
\\
\\
(c-0) Type $\Lambda_{12}^{0}$
\begin{align}
\label{o}
I(\Lambda_{12}^{0})&=\int_{X^{3}}\chi(x_{1})\chi(x_{2})\chi(x_{3})e^{(w_{1}x_{1}+w_{2}x_{2}+w_{3}x_{3}+w_{2}w_{3}y+w_{1}w_{3}y+w_{1}w_{2}y)t}dx_{1}dx_{2}dx_{3}
\\
\label{p}
&=\frac{w_{1}w_{2}w_{3}t^{3}e^{(w_{2}w_{3}+w_{1}w_{3}+w_{1}w_{2})yt}}{(e^{dw_{1}t}-1)(e^{dw_{2}t}-1)(e^{dw_{3}t}-1)}
\end{align}
\begin{align*}
\qquad \qquad \qquad \qquad\qquad \times (\sum_{a=0}^{d-1}\chi(a)e^{aw_{1}t})(\sum_{a=0}^{d-1}\chi(a)e^{aw_{2}t})(\sum_{a=0}^{d-1}\chi(a)e^{aw_{3}t})
\end{align*}
\\
\\
(c-1) Type $\Lambda_{12}^{1}$
\begin{align}
\label{q}
I(\Lambda_{12}^{1})&=\frac{d^{3}\int_{X^{3}}\chi(x_{1})\chi(x_{2})\chi(x_{3})e^{(w_{1}x_{1}+w_{2}x_{2}+w_{3}x_{3})t}dx_{1}dx_{2}dx_{3}}{\int_{X^{3}}e^{d(w_{2}w_{3}z_{1}+w_{1}w_{3}z_{2}+w_{1}w_{2}z_{3})t}dz_{1}dz_{2}dz_{3}}\qquad\qquad\qquad\qquad
\\
\label{r}
&=\frac{(e^{dw_{2}w_{3}t}-1)(e^{dw_{1}w_{3}t}-1)(e^{dw_{1}w_{2}t}-1)}{w_{1}w_{2}w_{3}(e^{dw_{1}t}-1)(e^{dw_{2}t}-1)(e^{dw_{3}t}-1)}
\end{align}
\begin{align*}
\qquad \qquad \qquad \qquad\qquad \times (\sum_{a=0}^{d-1}\chi(a)e^{aw_{1}t})(\sum_{a=0}^{d-1}\chi(a)e^{aw_{2}t})(\sum_{a=0}^{d-1}\chi(a)e^{aw_{3}t})
\end{align*}
All of the above $p$-adic integrals of various types are invariant
under all permutations of $w_{1},w_{2},w_{3}$, as one can see either
from $p$-adic integral representations in (\ref{k}), (\ref{m}),
(\ref{o}), and (\ref{q}) or from their explicit evaluations in
(\ref{l}), (\ref{n}), (\ref{p}), and (\ref{r}).

\section{Identities for generalized Bernoulli polynomials}
All of the following results can be easily obtained from (\ref{e})
and (\ref{h})-(\ref{j}). First, let's consider Type
$\Lambda_{23}^{i}$, for each $i=0,1,2,3.$
\\
\\
(a-0)
\begin{equation*}
\begin{split}
I(\Lambda_{23}^{0})=\int_{X}\chi(x_{1})e^{w_{2}w_{3}(x_{1}+w_{1}y_{1})t}dx_{1}\int_{X}\chi(x_{2})&e^{w_{1}w_{3}(x_{2}+w_{2}y_{2})t}dx_{2}\qquad\qquad\qquad\qquad\\
&\times \int_{X}\chi(x_{3})e^{w_{1}w_{2}(x_{3}+w_{3}y_{3})t}dx_{3}
\end{split}
\end{equation*}
\begin{equation*}
=(\sum_{k=0}^{\infty}\frac{B_{k,\chi}(w_{1}y_{1})}{k!}(w_{2}w_{3}t)^{k})(\sum_{l=0}^{\infty}\frac{B_{l,\chi}(w_{2}y_{2})}{l!}(w_{1}w_{3}t)^{l})(\sum_{m=0}^{\infty}\frac{B_{m,\chi}(w_{3}y_{3})}{m!}(w_{1}w_{2}t)^{m})\\
\end{equation*}
\begin{equation}\label{s}
=\sum_{n=0}^{\infty}(\sum_{k+l+m=n}^{}\binom{n}{k,l,m}B_{k,\chi}(w_{1}y_{1})B_{l,\chi}(w_{2}y_{2})B_{m,\chi}(w_{3}y_{3})w_{1}^{l+m}w_{2}^{k+m}w_{3}^{k+l})\frac{t^{n}}{n!},\qquad\qquad\\
\end{equation}
where the inner sum is over all nonnegative integers $k,l,m$, with
$k+l+m=n$, and
\begin{equation*}
\binom{n}{k,l,m}=\frac{n!}{k!l!m!}.
\end{equation*}
\\
\\
(a-1) Here we write $I(\Lambda_{23}^{1})$ in two different ways:
\\
\\
(1)
\begin{equation}\label{t}
\begin{split}
I(\Lambda_{23}^{1})=\frac{1}{w_{3}}\int_{X}\chi(x_{1})e^{w_{2}w_{3}(x_{1}+w_{1}y_{1})t}dx_{1}&\int_{X}\chi(x_{2})e^{w_{1}w_{3}(x_{2}+w_{2}y_{2})t}dx_{2}\qquad\qquad\qquad\qquad\\
&\times
\frac{dw_{3}\int_{X}\chi(x_{3})e^{w_{1}w_{2}x_{3}t}dx_{3}}{\int_{X}e^{dw_{1}w_{2}w_{3}x_{4}t}dx_{4}}
\end{split}
\end{equation}
\begin{equation*}
\begin{split}
=\frac{1}{w_{3}}(\sum_{k=0}^{\infty}B_{k,\chi}(w_{1}y_{1})\frac{(w_{2}w_{3}t)^{k}}{k!})(\sum_{l=0}^{\infty}&B_{l,\chi}(w_{2}y_{2})\frac{(w_{1}w_{3}t)^{l}}{l!})\qquad\qquad\\
&\times(\sum_{m=0}^{\infty}S_{m}(w_{3}d-1,\chi)\frac{(w_{1}w_{2}t)^{m}}{m!})
\end{split}
\end{equation*}
\begin{equation}\label{u}
=\sum_{n=0}^{\infty}(\sum_{k+l+m=n}^{}\binom{n}{k,l,m}B_{k,\chi}(w_{1}y_{1})B_{l,\chi}(w_{2}y_{2})S_{m}(w_{3}d-1,\chi)w_{1}^{l+m}w_{2}^{k+m}w_{3}^{k+l-1})\frac{t^{n}}{n!}.
\end{equation}
\\
(2) Invoking (\ref{i}), (\ref{t}) can also be written as\\

$I(\Lambda_{23}^{1})$
\begin{equation*}
=\frac{1}{w_{3}}\sum_{a=0}^{w_{3}d-1}\chi(a)\int_{X}\chi(x_{1})e^{w_{2}w_{3}(x_{1}+w_{1}y_{1})t}dx_{1}\int_{X}\chi(x_{2})e^{w_{1}w_{3}(x_{2}+w_{2}y_{2}+\frac{w_{2}}{w_{3}}a)t}dx_{2}\\
\end{equation*}
\begin{equation*}
=\frac{1}{w_{3}}\sum_{a=0}^{w_{3}d-1}\chi(a)(\sum_{k=0}^{\infty}B_{k,\chi}(w_{1}y_{1})\frac{(w_{2}w_{3}t)^{k}}{k!})(\sum_{l=0}^{\infty}B_{l,\chi}(w_{2}y_{2}+\frac{w_{2}}{w_{3}}a)\frac{(w_{1}w_{3}t)^{l}}{l!})\quad\qquad\\
\end{equation*}
\begin{equation}\label{v}
=\sum_{n=0}^{\infty}(w_{3}^{n-1}\sum_{k=0}^{n}\binom{n}{k}B_{k,\chi}(w_{1}y_{1})\sum_{a=0}^{w_{3}d-1}\chi(a)B_{n-k,\chi}(w_{2}y_{2}+\frac{w_{2}}{w_{3}}a)w_{1}^{n-k}w_{2}^{k})\frac{t^{n}}{n!}.\quad
\end{equation}
\\
(a-2) Here we write $I(\Lambda_{23}^{2})$ in three different ways:
\\
(1)
\begin{equation}\label{w}
I(\Lambda_{23}^{2})=\frac{1}{w_{2}w_{3}}\int_{X}\chi(x_{1})e^{w_{2}w_{3}(x_{1}+w_{1}y_{1})t}dx_{1}
\times
\frac{dw_{2}\int_{X}\chi(x_{2})e^{w_{1}w_{3}x_{2}t}dx_{2}}{\int_{X}e^{dw_{1}w_{2}w_{3}x_{4}t}dx_{4}}\qquad\qquad\qquad\qquad\qquad
\end{equation}
\begin{equation*}
\qquad\qquad\qquad\qquad\qquad\qquad\qquad\qquad\qquad\qquad\qquad\times
\frac{dw_{3}\int_{X}\chi(x_{3})e^{w_{1}w_{2}x_{3}t}dx_{3}}{\int_{X}e^{dw_{1}w_{2}w_{3}x_{4}t}dx_{4}}
\end{equation*}
\begin{equation*}
=\frac{1}{w_{2}w_{3}}(\sum_{k=0}^{\infty}{B_{k,\chi}(w_{1}y_{1})\frac{(w_{2}w_{3}t)^k}{k!}})
(\sum_{l=0}^{\infty}{S_{l}(w_{2}d-1,\chi)\frac{(w_{1}w_{3}t)^l}{l!}})\qquad
\end{equation*}
\begin{equation*}
\qquad\qquad\qquad\qquad\qquad\qquad\qquad\qquad\qquad\qquad\times(\sum_{m=0}^{\infty}{S_{m}(w_{3}d-1,\chi)\frac{(w_{1}w_{2}t)^m}{m!}})\quad\quad
\end{equation*}
\begin{equation}\label{x}
\begin{split}
=\sum_{n=0}^{\infty}(\sum_{k+l+m=n}\binom{n}{k,l,m}&B_{k,\chi}(w_{1}y_{1})S_{l}(w_{2}d-1,\chi)\qquad\qquad\qquad\qquad\qquad \\
&\times
S_{m}(w_{3}d-1,\chi)w_{1}^{l+m}w_{2}^{k+m-1}w_{3}^{k+l-1})\frac{t^{n}}{n!}.
\end{split}
\end{equation}
\\
\\
(2) Invoking (\ref{i}), (\ref{w}) can also be written as
\\
\\
~$I(\Lambda_{23}^{2})$
\begin{equation}\label{y}
=\frac{1}{w_{2}w_{3}}\sum_{a=0}^{w_{2}d-1}\chi(a)\int_{X}\chi(x_{1})e^{w_{2}w_{3}(x_{1}+w_{1}y_{1}+\frac{w_{1}}{w_{2}}a)t}dx_{1}
\times
\frac{dw_{3}\int_{X}\chi(x_{3})e^{w_{1}w_{2}x_{3}t}dx_{3}}{\int_{X}e^{dw_{1}w_{2}w_{3}x_{4}t}dx_{4}}
\end{equation}
\begin{equation*}
=\frac{1}{w_{2}w_{3}}\sum_{a=0}^{w_{2}d-1}\chi(a)(\sum_{k=0}^{\infty}B_{k,\chi}(w_{1}y_{1}+\frac{w_{1}}{w_{2}}a)\frac{(w_{2}w_{3}t)^k}{k!})
(\sum_{l=0}^{\infty}S_{l}(w_{3}d-1,\chi)\frac{(w_{1}w_{2}t)^l}{l!})
\end{equation*}
\begin{equation}\label{z}
=\sum_{n=0}^{\infty}(w_{2}^{n-1}\sum_{k=0}^{n}\binom{n}{k}\sum_{a=0}^{w_{2}d-1}\chi(a)B_{k,\chi}(w_{1}y_{1}+\frac{w_{1}}{w_{2}}a)S_{n-k}(w_{3}d-1,\chi)w_{1}^{n-k}w_{3}^{k-1})\frac{t^{n}}{n!}.\quad
\end{equation}
\\
\\
(3) Invoking (\ref{i}) once again, (\ref{y}) can be written as
\\
\\
\begin{equation*}
I(\Lambda_{23}^{2})=\frac{1}{w_{2}w_{3}}\sum_{a=0}^{w_{2}d-1}\chi(a)\sum_{b=0}^{w_{3}d-1}\chi(b)\int_{X}\chi(x_{1})e^{w_{2}w_{3}(x_{1}+w_{1}y_{1}+\frac{w_{1}}{w_{2}}a+\frac{w_{1}}{w_{3}}b)t}dx_{1}~
\end{equation*}
\begin{equation*}
\qquad\qquad=\frac{1}{w_{2}w_{3}}\sum_{a=0}^{w_{2}d-1}\chi(a)\sum_{b=0}^{w_{3}d-1}\chi(b)(\sum_{n=0}^{\infty}B_{n,\chi}(w_{1}y_{1}+\frac{w_{1}}{w_{2}}a+\frac{w_{1}}{w_{3}}b)\frac{(w_{2}w_{3}t)^n}{n!})
\end{equation*}
\begin{equation}\label{a1}
=\sum_{n=0}^{\infty}((w_{2}w_{3})^{n-1}\sum_{a=0}^{w_{2}d-1}\sum_{b=0}^{w_{3}d-1}\chi(ab)B_{n,\chi}(w_{1}y_{1}+\frac{w_{1}}{w_{2}}a+\frac{w_{1}}{w_{3}}b))\frac{t^{n}}{n!}.\quad
\end{equation}
\\
\\
(a-3)
\\
\\
\begin{equation*}
I(\Lambda_{23}^{3})=\frac{1}{w_{1}w_{2}w_{3}}\times
\frac{dw_{1}\int_{X}\chi(x_{1})e^{w_{2}w_{3}x_{1}t}dx_{1}}{\int_{X}e^{dw_{1}w_{2}w_{3}x_{4}t}dx_{4}}
\times\frac{dw_{2}\int_{X}\chi(x_{2})e^{w_{1}w_{3}x_{2}t}dx_{2}}{\int_{X}e^{dw_{1}w_{2}w_{3}x_{4}t}dx_{4}}\qquad\qquad
\end{equation*}
\begin{equation*}
\qquad\qquad\qquad\qquad\times\frac{dw_{3}\int_{X}\chi(x_{3})e^{w_{1}w_{2}x_{3}t}dx_{3}}{\int_{X}e^{dw_{1}w_{2}w_{3}x_{4}t}dx_{4}}
\end{equation*}
\begin{align*}
\quad\qquad=\frac{1}{w_{1}w_{2}w_{3}}(\sum_{k=0}^{\infty}S_{k}(w_{1}d-1,\chi)\frac{(w_{2}w_{3}t)^{k}}{k!})
&(\sum_{l=0}^{\infty}S_{l}(w_{2}d-1,\chi)\frac{(w_{1}w_{3}t)^{l}}{l!})\quad\qquad\\
&\times(\sum_{m=0}^{\infty}S_{m}(w_{3}d-1,\chi)\frac{(w_{1}w_{2}t)^{m}}{m!})
\end{align*}
\begin{equation}\label{b1}
\begin{split}
=\sum_{n=0}^{\infty}\sum_{k+l+m=n}^{}(\binom{n}{k,l,m}S_{k}(w_{1}d-1,\chi)S_{l}(w_{2}d-1,&\chi)S_{m}(w_{3}d-1,\chi)\\
&\times w_{1}^{l+m-1}w_{2}^{k+m-1}w_{3}^{k+l-1})\frac{t^{n}}{n!}.\quad
\end{split}
\end{equation}
\\
(b) For Type $\Lambda_{13}^{i}~(i=0,1,2,3)$, we may consider the analogous things to the
ones in (a-0), (a-1), (a-2), and (a-3). However, these do not lead us to new identities.
Indeed, if we substitute $w_{2}w_{3},w_{1}w_{3},w_{1}w_{2}$ respectively for $w_{1},w_{2},w_{3}$
in (\ref{k}), this amounts to replacing $t$ by $w_{1}w_{2}w_{3}t$ in (\ref{m}).
So, upon replacing  $w_{1},w_{2},w_{3}$ respectively by $w_{2}w_{3},w_{1}w_{3},w_{1}w_{2}$
and dividing by $(w_{1}w_{2}w_{3})^n$, in each of the expressions
of (\ref{s}), (\ref{u}), (\ref{v}), (\ref{x}), (\ref{z})-(\ref{b1}), we will get the corresponding symmetric identities
for  Type $\Lambda_{13}^{i}~(i=0,1,2,3)$.
\\
\\
(c-0)
\begin{align*}
&I(\Lambda_{12}^{0})\\
&=\int_{X}\chi(x_{1})e^{w_{1}(x_{1}+w_{2}y)t}dx_{1}
\int_{X}\chi(x_{2})e^{w_{2}(x_{2}+w_{3}y)t}dx_{2}\int_{X}\chi(x_{3})e^{w_{3}(x_{3}+w_{1}y)t}dx_{3}\\
&=(\sum_{n=0}^{\infty}\frac{B_{k,\chi}(w_{2}y)}{k!}(w_{1}t)^{k})(\sum_{l=0}^{\infty}\frac{B_{l,\chi}(w_{3}y)}{l!}(w_{2}t)^{l})(\sum_{m=0}^{\infty}\frac{B_{m,\chi}(w_{1}y)}{m!}(w_{3}t)^{m})
\end{align*}
\begin{align}\label{c1}
=\sum_{n=0}^{\infty}(\sum_{k+l+m=n}\binom{n}{k,l,m}B_{k,\chi}(w_{2}y)B_{l,\chi}(w_{3}y)B_{m,\chi}(w_{1}y)w_{1}^{k}w_{2}^{l}w_{3}^{m})\frac{t^n}{n!}.~\qquad\qquad
\end{align}
\\
\\
(c-1)
\begin{align*}
I(\Lambda_{12}^{1})=\frac{1}{w_{1}w_{2}w_{3}}\frac{dw_{2}\int_{X}\chi(x_{1})e^{w_{1}x_{1}t}dx_{1}}{\int_{X}e^{dw_{1}w_{2}z_{3}t}dz_{3}}
\times\frac{dw_{3}\int_{X}\chi(x_{2})e^{w_{2}x_{2}t}dx_{2}}{\int_{X}e^{dw_{2}w_{3}z_{1}t}dz_{1}}\qquad\qquad\\
\times\frac{dw_{1}\int_{X}\chi(x_{3})e^{w_{3}x_{3}t}dx_{3}}{\int_{X}e^{dw_{3}w_{1}z_{2}t}dz_{2}}\\
=\frac{1}{w_{1}w_{2}w_{3}}(\sum_{k=0}^{\infty}S_{k}(w_{2}d-1,\chi)\frac{(w_{1}t)^{k}}{k!})
(\sum_{l=0}^{\infty}S_{l}(w_{3}d-1,\chi)\frac{(w_{2}t)^{l}}{l!})\qquad\\
\times
(\sum_{m=0}^{\infty}S_{m}(w_{1}d-1,\chi)\frac{(w_{3}t)^{m}}{m!})
\end{align*}
\begin{equation}\label{d1}
\begin{split}
=\sum_{n=0}^{\infty}(\sum_{k+l+m=n}^{}\binom{n}{k,l,m}S_{k}(w_{2}d-1,\chi)S_{l}(w_{3}d-1,\chi)S_{m}&(w_{1}d-1,\chi)\qquad\qquad\\
&\times w_{1}^{k-1}w_{2}^{l-1}w_{3}^{m-1})\frac{t^{n}}{n!}.
\end{split}
\end{equation}
\section{Main theorems}
As we noted earlier in the last paragraph of Section 2, the various
types of quotients of $p$-adic integrals are invariant under any
permutation of $w_{1},w_{2},w_{3}$. So the corresponding expressions
in Section 3 are also invariant under any permutation of
$w_{1},w_{2},w_{3}$. Thus our results about identities of symmetry
will be immediate consequences of this observation.

However, not all permutations of an expression in Section 3 yield
distinct ones. In fact, as these expressions are obtained by
permuting $w_{1},w_{2},w_{3}$ in a single one labeled by them, they
can be viewed as a group in a natural manner and hence it is
isomorphic to a quotient of $S_{3}$. In particular, the number of
possible distinct expressions are $1,2,3$ or $6$. (a-0), (a-1(1)),
(a-1(2)), and (a-2(2)) give the full six identities of symmetry,
(a-2(1)) and (a-2(3)) yield three identities of symmetry, and (c-0)
and (c-1) give two identities of symmetry, while the expression in
(a-3) yields no identities of symmetry.

Here we will just consider the cases of Theorems 4 and 8, leaving the others as easy
exercises for the reader. As for the case of Theorem 4, in addition to (\ref{o1})-(\ref{q1}),
we get the following three ones:\\
\begin{align}
\label{e1}
&\sum_{k+l+m=n}\binom{n}{k,l,m}B_{k,\chi}(w_{1}y_{1})S_{l}(w_{3}d-1,\chi)S_{m}(w_{2}d-1,\chi)w_{1}^{l+m}w_{3}^{k+m-1}w_{2}^{k+l-1},\\
\label{f1}
&\sum_{k+l+m=n}\binom{n}{k,l,m}B_{k,\chi}(w_{2}y_{1})S_{l}(w_{1}d-1,\chi)S_{m}(w_{3}d-1,\chi)w_{2}^{l+m}w_{1}^{k+m-1}w_{3}^{k+l-1},\\
\label{g1}
&\sum_{k+l+m=n}\binom{n}{k,l,m}B_{k,\chi}(w_{3}y_{1})S_{l}(w_{2}d-1,\chi)S_{m}(w_{1}d-1,\chi)w_{3}^{l+m}w_{2}^{k+m-1}w_{1}^{k+l-1}.
\end{align}
But, by interchanging $l$ and $m$, we see that (\ref{e1}),
(\ref{f1}), and (\ref{g1}) are respectively equal to (\ref{o1}),
(\ref{p1}), and (\ref{q1}).
\\
As to Theorem 8, in addition to (\ref{u1}) and (\ref{v1}), we have:
\begin{align}
\label{h1}
&\sum_{k+l+m=n}\binom{n}{k,l,m}S_{k}(w_{2}d-1,\chi)S_{l}(w_{3}d-1,\chi)S_{m}(w_{1}d-1,\chi)w_{1}^{k-1}w_{2}^{l-1}w_{3}^{m-1},\\
\label{i1}
&\sum_{k+l+m=n}\binom{n}{k,l,m}S_{k}(w_{3}d-1,\chi)S_{l}(w_{1}d-1,\chi)S_{m}(w_{2}d-1,\chi)w_{2}^{k-1}w_{3}^{l-1}w_{1}^{m-1},\\
\label{j1}
&\sum_{k+l+m=n}\binom{n}{k,l,m}S_{k}(w_{3}d-1,\chi)S_{l}(w_{2}d-1,\chi)S_{m}(w_{1}d-1,\chi)w_{1}^{k-1}w_{3}^{l-1}w_{2}^{m-1},\\
\label{k1}
&\sum_{k+l+m=n}\binom{n}{k,l,m}S_{k}(w_{2}d-1,\chi)S_{l}(w_{1}d-1,\chi)S_{m}(w_{3}d-1,\chi)w_{3}^{k-1}w_{2}^{l-1}w_{1}^{m-1}.
\end{align}
\\
  However, (\ref{h1}) and (\ref{i1}) are equal to (\ref{u1}), as we
can see by applying the permutations $k\rightarrow l,l\rightarrow
m,m\rightarrow k$ for (\ref{h1}) and $k\rightarrow m,l\rightarrow
k,m\rightarrow l$ for (\ref{i1}). Similarly, we see that (\ref{j1})
and (\ref{k1}) are equal to (\ref{v1}), by applying permutations
$k\rightarrow l,l\rightarrow m,m\rightarrow k$ for (\ref{j1}) and
$k\rightarrow m,l\rightarrow k,m\rightarrow l$ for (\ref{k1}).
\begin{theorem}\label{A}
Let $w_{1},w_{2},w_{3}$ be any positive integers. Then we have:
\begin{equation}\label{l1}
\begin{split}
&\sum_{k+l+m=n}\binom{n}{k,l,m}B_{k,\chi}(w_{1}y_{1})B_{l,\chi}(w_{2}y_{2})B_{m,\chi}(w_{3}y_{3})w_{1}^{l+m}w_{2}^{k+m}w_{3}^{k+l}\\
=&\sum_{k+l+m=n}\binom{n}{k,l,m}B_{k,\chi}(w_{1}y_{1})B_{l,\chi}(w_{3}y_{2})B_{m,\chi}(w_{2}y_{3})w_{1}^{l+m}w_{3}^{k+m}w_{2}^{k+l}\\
=&\sum_{k+l+m=n}\binom{n}{k,l,m}B_{k,\chi}(w_{2}y_{1})B_{l,\chi}(w_{1}y_{2})B_{m,\chi}(w_{3}y_{3})w_{2}^{l+m}w_{1}^{k+m}w_{3}^{k+l}\\
=&\sum_{k+l+m=n}\binom{n}{k,l,m}B_{k,\chi}(w_{2}y_{1})B_{l,\chi}(w_{3}y_{2})B_{m,\chi}(w_{1}y_{3})w_{2}^{l+m}w_{3}^{k+m}w_{1}^{k+l}\\
=&\sum_{k+l+m=n}\binom{n}{k,l,m}B_{k,\chi}(w_{3}y_{1})B_{l,\chi}(w_{1}y_{2})B_{m,\chi}(w_{2}y_{3})w_{3}^{l+m}w_{1}^{k+m}w_{2}^{k+l}\\
=&\sum_{k+l+m=n}\binom{n}{k,l,m}B_{k,\chi}(w_{3}y_{1})B_{l,\chi}(w_{2}y_{2})B_{m,\chi}(w_{1}y_{3})w_{3}^{l+m}w_{2}^{k+m}w_{1}^{k+l}.
\end{split}
\end{equation}
\end{theorem}
\begin{theorem}\label{B}
  Let $w_{1},w_{2},w_{3}$ be any positive integers. Then we have:
\begin{equation}\label{m1}
\begin{split}
&\sum_{k+l+m=n}\binom{n}{k,l,m}B_{k,\chi}(w_{1}y_{1})B_{l,\chi}(w_{2}y_{2})S_{m}(w_{3}d-1,\chi)w_{1}^{l+m}w_{2}^{k+m}w_{3}^{k+l-1}\\
=&\sum_{k+l+m=n}\binom{n}{k,l,m}B_{k,\chi}(w_{1}y_{1})B_{l,\chi}(w_{3}y_{2})S_{m}(w_{2}d-1,\chi)w_{1}^{l+m}w_{3}^{k+m}w_{2}^{k+l-1}\\
=&\sum_{k+l+m=n}\binom{n}{k,l,m}B_{k,\chi}(w_{2}y_{1})B_{l,\chi}(w_{1}y_{2})S_{m}(w_{3}d-1,\chi)w_{2}^{l+m}w_{1}^{k+m}w_{3}^{k+l-1}\\
=&\sum_{k+l+m=n}\binom{n}{k,l,m}B_{k,\chi}(w_{2}y_{1})B_{l,\chi}(w_{3}y_{2})S_{m}(w_{1}d-1,\chi)w_{2}^{l+m}w_{3}^{k+m}w_{1}^{k+l-1}\\
=&\sum_{k+l+m=n}\binom{n}{k,l,m}B_{k,\chi}(w_{3}y_{1})B_{l,\chi}(w_{2}y_{2})S_{m}(w_{1}d-1,\chi)w_{3}^{l+m}w_{2}^{k+m}w_{1}^{k+l-1}\\
=&\sum_{k+l+m=n}\binom{n}{k,l,m}B_{k,\chi}(w_{3}y_{1})B_{l,\chi}(w_{1}y_{2})S_{m}(w_{2}d-1,\chi)w_{3}^{l+m}w_{1}^{k+m}w_{2}^{k+l-1}.
\end{split}
\end{equation}
\end{theorem}

\begin{theorem}\label{C}
Let $w_{1},w_{2},w_{3}$ be any positive integers. Then we have:
\begin{equation*}
\begin{split}
&w_{1}^{n-1}\sum_{k=0}^{n}\binom{n}{k}B_{k,\chi}(w_{3}y_{1})\sum_{a=0}^{w_{1}d-1}\chi(a)B_{n-k,\chi}(w_{2}y_{2}+\frac{w_{2}}{w_{1}}a)w_{3}^{n-k}w_{2}^{k}\\
=&w_{1}^{n-1}\sum_{k=0}^{n}\binom{n}{k}B_{k,\chi}(w_{2}y_{1})\sum_{a=0}^{w_{1}d-1}\chi(a)B_{n-k,\chi}(w_{3}y_{2}+\frac{w_{3}}{w_{1}}a)w_{2}^{n-k}w_{3}^{k}
\end{split}
\end{equation*}
\begin{equation}\label{n1}
\begin{split}
=&w_{2}^{n-1}\sum_{k=0}^{n}\binom{n}{k}B_{k,\chi}(w_{3}y_{1})\sum_{a=0}^{w_{2}d-1}\chi(a)B_{n-k,\chi}(w_{1}y_{2}+\frac{w_{1}}{w_{2}}a)w_{3}^{n-k}w_{1}^{k}\\
=&w_{2}^{n-1}\sum_{k=0}^{n}\binom{n}{k}B_{k,\chi}(w_{1}y_{1})\sum_{a=0}^{w_{2}d-1}\chi(a)B_{n-k,\chi}(w_{3}y_{2}+\frac{w_{3}}{w_{2}}a)w_{1}^{n-k}w_{3}^{k}\\
=&w_{3}^{n-1}\sum_{k=0}^{n}\binom{n}{k}B_{k,\chi}(w_{2}y_{1})\sum_{a=0}^{w_{3}d-1}\chi(a)B_{n-k,\chi}(w_{1}y_{2}+\frac{w_{1}}{w_{2}}a)w_{2}^{n-k}w_{1}^{k}\\
=&w_{3}^{n-1}\sum_{k=0}^{n}\binom{n}{k}B_{k,\chi}(w_{1}y_{1})\sum_{a=0}^{w_{3}d-1}\chi(a)B_{n-k,\chi}(w_{2}y_{2}+\frac{w_{2}}{w_{3}}a)w_{1}^{n-k}w_{2}^{k}.
\end{split}
\end{equation}
\end{theorem}

\begin{theorem}\label{D}
Let $w_{1},w_{2},w_{3}$ be any positive integers. Then we have the
following three symmetries in $w_{1},w_{2},w_{3}$:

\begin{equation}
\label{o1}
\sum_{k+l+m=n}\binom{n}{k,l,m}B_{k,\chi}(w_{1}y_{1})S_{l}(w_{2}d-1,\chi)S_{m}(w_{3}d-1,\chi)w_{1}^{l+m}w_{2}^{k+m-1}w_{3}^{k+l-1}
\end{equation}
\begin{equation}
\label{p1}
=\sum_{k+l+m=n}\binom{n}{k,l,m}B_{k,\chi}(w_{2}y_{1})S_{l}(w_{3}d-1,\chi)S_{m}(w_{1}d-1,\chi)w_{2}^{l+m}w_{3}^{k+m-1}w_{1}^{k+l-1}
\end{equation}
\begin{equation}
\label{q1}
=\sum_{k+l+m=n}\binom{n}{k,l,m}B_{k,\chi}(w_{3}y_{1})S_{l}(w_{1}d-1,\chi)S_{m}(w_{2}d-1,\chi)w_{3}^{l+m}w_{1}^{k+m-1}w_{2}^{k+l-1}.
\end{equation}
\end{theorem}

\begin{theorem}\label{E}
Let $w_{1},w_{2},w_{3}$ be any positive integers. Then we have:
\begin{equation}\label{r1}
\begin{split}
&w_{1}^{n-1}\sum_{k=0}^{n}\binom{n}{k}\sum_{a=0}^{w_{1}d-1}\chi(a)B_{k,\chi}(w_{2}y_{1}+\frac{w_{2}}{w_{1}}a)S_{n-k}(w_{3}d-1,\chi)w_{2}^{n-k}w_{3}^{k-1}\\
=&w_{1}^{n-1}\sum_{k=0}^{n}\binom{n}{k}\sum_{a=0}^{w_{1}d-1}\chi(a)B_{k,\chi}(w_{3}y_{1}+\frac{w_{3}}{w_{1}}a)S_{n-k}(w_{2}d-1,\chi)w_{3}^{n-k}w_{2}^{k-1}\\
=&w_{2}^{n-1}\sum_{k=0}^{n}\binom{n}{k}\sum_{a=0}^{w_{2}d-1}\chi(a)B_{k,\chi}(w_{1}y_{1}+\frac{w_{1}}{w_{2}}a)S_{n-k}(w_{3}d-1,\chi)w_{1}^{n-k}w_{3}^{k-1}\\
=&w_{2}^{n-1}\sum_{k=0}^{n}\binom{n}{k}\sum_{a=0}^{w_{2}d-1}\chi(a)B_{k,\chi}(w_{3}y_{1}+\frac{w_{3}}{w_{2}}a)S_{n-k}(w_{1}d-1,\chi)w_{3}^{n-k}w_{1}^{k-1}\\
=&w_{3}^{n-1}\sum_{k=0}^{n}\binom{n}{k}\sum_{a=0}^{w_{3}d-1}\chi(a)B_{k,\chi}(w_{1}y_{1}+\frac{w_{1}}{w_{3}}a)S_{n-k}(w_{2}d-1,\chi)w_{1}^{n-k}w_{2}^{k-1}\\
=&w_{3}^{n-1}\sum_{k=0}^{n}\binom{n}{k}\sum_{a=0}^{w_{3}d-1}\chi(a)B_{k,\chi}(w_{2}y_{1}+\frac{w_{2}}{w_{3}}a)S_{n-k}(w_{1}d-1,\chi)w_{2}^{n-k}w_{1}^{k-1}.
\end{split}
\end{equation}
\end{theorem}

\begin{theorem}\label{F}
Let $w_{1},w_{2},w_{3}$ be any positive integers. Then we have the
following three symmetries in $w_{1},w_{2},w_{3}$:\\
\begin{equation}\label{s1}
\begin{split}
&(w_{1}w_{2})^{n-1}\sum_{a=0}^{w_{1}d-1}\sum_{b=0}^{w_{2}d-1}\chi(ab)B_{n,\chi}(w_{3}y_{1}+\frac{w_{3}}{w_{1}}a+\frac{w_{3}}{w_{2}}b)\\
=&(w_{2}w_{3})^{n-1}\sum_{a=0}^{w_{2}d-1}\sum_{b=0}^{w_{3}d-1}\chi(ab)B_{n,\chi}(w_{1}y_{1}+\frac{w_{1}}{w_{2}}a+\frac{w_{1}}{w_{3}}b)\\
=&(w_{3}w_{1})^{n-1}\sum_{a=0}^{w_{3}d-1}\sum_{b=0}^{w_{1}d-1}\chi(ab)B_{n,\chi}(w_{2}y_{1}+\frac{w_{2}}{w_{3}}a+\frac{w_{2}}{w_{1}}b).
\end{split}
\end{equation}
\end{theorem}

\begin{theorem}\label{G}
Let $w_{1},w_{2},w_{3}$ be any positive integers. Then we have the
following two symmetries in $w_{1},w_{2},w_{3}$:\\
\begin{equation}\label{t1}
\begin{split}
&\sum_{k+l+m=n}\binom{n}{k,l,m}B_{k,\chi}(w_{1}y)B_{l,\chi}(w_{2}y)B_{m,\chi}(w_{3}y)w_{3}^{k}w_{1}^{l}w_{2}^{m}\\
=&\sum_{k+l+m=n}\binom{n}{k,l,m}B_{k,\chi}(w_{1}y)B_{l,\chi}(w_{3}y)B_{m,\chi}(w_{2}y)w_{2}^{k}w_{1}^{l}w_{3}^{m}.
\end{split}
\end{equation}
\end{theorem}

\begin{theorem}\label{H}
Let $w_{1},w_{2},w_{3}$ be any positive integers. Then we have
the following two symmetries in $w_{1},w_{2},w_{3}$:
\begin{align}
\label{u1}
&\sum_{k+l+m=n}\binom{n}{k,l,m}S_{k}(w_{1}d-1,\chi)S_{l}(w_{2}d-1,\chi)S_{m}(w_{3}d-1,\chi)w_{3}^{k-1}w_{1}^{l-1}w_{2}^{m-1}\\
\label{v1}
=&\sum_{k+l+m=n}\binom{n}{k,l,m}S_{k}(w_{1}d-1,\chi)S_{l}(w_{3}d-1,\chi)S_{m}(w_{2}d-1,\chi)w_{2}^{k-1}w_{1}^{l-1}w_{3}^{m-1}.
\end{align}
\end{theorem}

\bibliographystyle{amsplain}

\end{document}